\documentclass[a4paper,10pt]{amsart}

\usepackage{amsmath}
\usepackage{amsthm}
\usepackage{amsfonts}
\usepackage{amssymb}
\usepackage[all]{xy}
\usepackage{nicefrac}

\newcommand{\myAA}{\mathcal{A}}
\newcommand{\AAA}{R\langle\langle A\rangle\rangle}
\newcommand{\myFF}{\mathcal{F}}

\newcommand{\myCC}{\mathcal{C}}

\newcommand{\ZZ}{\mathbb{Z}}
\def\CC{\mathbb{C}}

\newcommand{\PP}{\mathbb{P}}
\newcommand{\QQ}{\mathbb{Q}}

\newcommand{\BB}{\mathbf{B}}
\newcommand{\TT}{\mathbf{T}}

\newtheorem{theorem}{Theorem}[section]

\newtheorem{lemma}{Lemma}[section]
\newtheorem{remark}{Remark}[section]

\title{Positivity in skew--symmetric cluster algebras of finite type}
\author{Giovanni Cerulli Irelli}
\address{Sapienza Universit\`a di Roma. Piazzale Aldo Moro 5, 00185, Rome (ITALY)}
\email{cerulli.math@googlemail.com}
\begin{document}

\begin{abstract}
We prove that the basis of cluster monomials of a skew--symmetric cluster  algebra $\myAA$ of finite type is the atomic basis of $\myAA$.  This means that an element of $\myAA$ is positive if and only if it has a non--negative expansion in the basis of cluster monomials. In particular cluster monomials are positive indecomposable, i.e. they cannot be written as a sum of positive elements. 
\textbf{Math. Subj. Class.(2010)}: 13F60.
\end{abstract}
\maketitle

\section{Introduction}

Let $H$ be an orientation of a simply--laced Dynkin diagram, i.e. a diagram of type A, D, or E in the Cartan--Killing classification. We consider the coefficient--free cluster algebra $\myAA(H)$ associated with $H$. This is a $\ZZ$--subalgebra of the field $\myFF=\QQ(x_1,\dotsc, x_n)$, where $n$ is the number of vertices of $H$, introduced by Fomin and Zelevinsky \cite{FZI}. The algebra $\myAA(H)$ can be described as follows: for every $k\in[1,n]$ let us consider the  element $x_k'\in\myFF$ defined by:
\begin{equation}\label{eq:MutIntro}
x_k'=\frac{\prod_{k\rightarrow j\in H_1}x_j+\prod_{i\rightarrow k\in H_1}x_i}{ x_k}
\end{equation}
(here $H_1$ denotes the set of arrows of H).  We have $\myAA(H)=\ZZ[x_k,x_k':\,k=1,\dotsc,n]\subset\myFF$ (see \cite[theorem~1.18, corollary~1.19]{BFZIII} or the survey \cite[theorem~4.13]{FZNotes}). Moreover, in \cite{BFZIII} it is  shown that standard monomials, i.e. monomials in $x_1,x_1',\dotsc, x_n,x'_n$ which do not contain the products $x_kx'_k$ ($k\in[1,n]$), form a $\ZZ$--basis of $\myAA(H)$.

Besides the basis of standard monomials, there is another basis of $\myAA(H)$ which is of particular interest to us. This is the basis $\BB$ of cluster monomials. Let us briefly define the set $\BB$ (see section~\ref{Sec:CluAlg} for more details). By definition, the algebra $\myAA(H)$ is the $\ZZ$--subalgebra generated by some elements of $\myFF$ called cluster variables. The cluster variables are grouped into free--generating sets for the field $\myFF$ called clusters. In particular, every cluster $\mathbf{s}$ consists of $n$ algebraically independent rational functions $s_1,\dotsc, s_n$ and $\myFF\simeq \QQ(s_1,\dotsc, s_n)$. A cluster monomial of $\myAA(H)$ is, by definition, a monomial in cluster variables belonging to the same cluster.

Cluster monomials are natural elements to consider in the additive categorification of cluster algebras via cluster categories \cite{BMRRT}. Namely they correspond to cluster--tilting objects. This description allow Caldero and Keller in \cite{CK1} to prove that cluster monomials form a $\ZZ$--basis of $\myAA(H)$. Cluster monomials are also important elements to consider in some geometric realizations of cluster algebras \cite{BFZIII}, \cite{FZII}, where every cluster provides a criterion for total positivity. Moreover cluster monomials belong to the dual semi--canonical basis of the coordinate ring $\CC[N]$ of a maximal unipotent group N whose Lie algebra $\mathfrak{n}$ is the maximal unipotent subalgebra of a simple Lie algebra of type A, D or E  \cite[theorem~2.8]{GLS06}. In this paper we give further evidence of the importance of cluster monomials in the theory of cluster algebras itself, as conjectured by Fomin and Zelevinsky \cite[conjecture~4.19]{FZNotes}. As a special case of the Fomin--Zelevinsky's Laurent phenomenon \cite{FZI}, every element of $\myAA(H)$ is a Laurent polynomial in all its clusters. An element $p$ of the cluster algebra $\myAA(H)$ is called positive if its Laurent expansions in all the clusters of $\myAA(H)$ have non--negative integer coefficients. Positive elements form a semiring, i.e. sums and products of positive elements are positive. We say (see \cite{Sherman}, \cite{CanonicalBasis}) that a $\ZZ$--basis $\BB$ of $\myAA(H)$ is atomic if  the semiring of positive elements consists precisely of $\ZZ_{\geq0}$--linear combinations of elements of $\BB$. Note that if an atomic basis exists, it is unique and it is formed by positive indecomposable elements, i.e. those elements which cannot be written as a sum of positive elements.

\begin{theorem}\label{MainThmIntro}
The set of cluster monomials is the atomic basis of $\myAA(H)$ when $H$ is an orientation of a simply--laced  Dynkin diagram.
\end{theorem}

The proof of theorem~\ref{MainThmIntro} is an application of the theory of quivers with potentials developed by Derksen, Weyman and Zelevinsky in \cite{DWZ} and \cite{DWZII}. In particular, the proof does not depend on the choice of  ``coefficients": for every semifield $\PP$ a cluster algebra $\myAA_\PP(H)$ is defined, which is said to have coefficients in $\PP$. Theorem~\ref{MainThmIntro} holds in $\myAA_\PP(H)$ as long as the cluster monomials form a $\ZZ\PP$--basis of it  (see remark~\ref{Rem:Coeff}). It is only for simplifying  notation that we prefer to work in the coefficient--free setting.

The existence of an atomic basis of a general cluster algebra is still an open problem. Only a few cases are known \cite{Sherman}, \cite{Thesis}, \cite{CanonicalBasis}.  We trust that our proof of theorem~\ref{MainThmIntro}, together with the techniques developed in \cite{CEsp} and \cite{CeDEsp}, is an important step for the description of the positive semiring of a cluster algebra associated with any acyclic  quiver.

We prove theorem~\ref{MainThmIntro} in section~\ref{Sec:Proof}. The other sections are devoted to recalling the results used in the proof.

\textbf{Acknowledgments}
This paper was written while I was at the Department of Mathematics of the University ``La~Sapienza" of Rome as a postdoctoral fellow under the supervision of Corrado~De~Concini.

The last version of the paper was written while I was part of the trimester program ``On the Interaction of Representation Theory with Geometry and Combinatorics" of the Hausdorff Research Institute for Mathematics of Bonn.

I thank both these institutions for  financial support and for excellent working conditions.

\section{Background on cluster algebras}\label{Sec:CluAlg}
In this section we recall the definition of a cluster algebra and some properties from \cite{FZI}.

Let $n\geq1$ be a positive integer. We consider an $n$--regular tree $\mathbf{T}_n$ and we label its edges by numbers $1,2\dotsc, n$ so that two edges adjacent to the same vertex receive different labels. We introduce the dynamic of seed mutations on $\mathbf{T}_n$. Let $\myFF=\QQ(x_1,\dotsc, x_n)$ be the field of rational functions in $n$ independent variables. A seed in $\myFF$ is a pair $(B,\mathbf{u})$ where $B$ is a skew--symmetric $n\times n$ integer matrix and $\mathbf{u}=(u_1,\dotsc, u_n)$ is a free--generating set for $\myFF$ so that $\myFF\simeq \QQ(u_1,\dotsc, u_n)$. The matrix $B$ is called the exchange matrix of the seed $\Sigma$, while the set $\mathbf{u}$ is called its cluster. The elements of the cluster of $\Sigma$ are called its cluster variables.
Given a seed $\Sigma$ of $\myFF$ and $k\in [1,n]$ (as is customary, we use the notation $[1,n]:=(1,2,\dotsc,n)$) we define a new seed $\mu_k(\Sigma)=(B',\mathbf{u}')$, called the mutation of $\Sigma$ in direction $k$, obtained from $\Sigma$ by the following rules of mutation:
\begin{enumerate}
\item the matrix $B'=(b'_{ij})$ is given by
\begin{equation}\label{Eq:MatrixMut}
b'_{ij}=\left\{
\begin{array}{lc}
-b_{ij}&\textrm{if }i=k\textrm{ or }j=k;\\
b_{ij}+sg(b_{ik})[b_{ik}b_{kj}]_+&\textrm{otherwise}
\end{array}\right.
\end{equation}
where $sg(b)$ denotes the sign of the integer $b$ and $[b]_+:=\textrm{max}(b,0)$.
\item The new cluster $\mathbf{u}'$ is obtained from the cluster $\mathbf{u}=(u_1,\dotsc, u_n)$ by $\mathbf{u}'=\mathbf{u}\setminus\{u_k\}\cup\{u'_k\}$ where
\begin{equation}\label{Eq:MutationUk}
u'_k=\frac{\prod_{i=1}^n u_i^{[b_{ik}]_+}+\prod_{j=1}^n u_j^{[-b_{jk}]_+}}{u_k}.
\end{equation}
\end{enumerate}
A cluster pattern is the assignment of a seed $\Sigma_t$ to every vertex of $\TT_n$ so that whenever $\xymatrix@1{t\ar@{-}^{k}[r]&t'}$, i.e. the unique edge adjacent to $t$ and labelled with $k$ connects $t$ with the vertex $t'$,  the assigned seeds $\Sigma_t$ and $\Sigma_{t'}$ satisfy $\Sigma_{t'}=\mu_k(\Sigma_t)$. It is clear that a cluster pattern  is uniquely determined by the choice of an ``initial" seed $\Sigma_0$ and we  denote it by $\TT_n(\Sigma_0)$. By definition, the (coefficient--free skew--symmetric) cluster algebra $\myAA(\Sigma_0)=\myAA(\TT_n(\Sigma_0))$ is the $\ZZ$--subalgebra of $\myFF$ generated by the cluster variables of the seeds of $\TT_n(\Sigma_0)$.

We notice that the cluster pattern depends uniquely on the choice of the initial exchange matrix $B$ of $\Sigma_0$ and we hence often write $\myAA(B)$  instead of $\myAA(\Sigma_0)$ and $\TT_n(B)$ instead of $\TT_n(\Sigma_0)$.

We sometimes prefer to use the language of quivers instead of the one of matrices.
Let $Q=(Q_0,Q_1,t,h)$ be a finite quiver without loops and oriented $2$--cycles, with vertex set $Q_0$ of cardinality $n$, with edges $Q_1$ and orientation given by the two maps $t,h:Q_1\rightarrow Q_0$ which associate to an edge $a$ its tail $t(a)$ and its head $h(a)$ and we write $\xymatrix@1{t(a)\ar^a[r]&h(a)}$. We associate with $Q$ the skew--symmetric $n\times n$ integer matrix $B(Q)$ whose $ij$--th entry equals the number of arrows from the vertex $j$ to the vertex $i$ minus the number of arrows from $i$ to $j$. The map $Q\mapsto B(Q)$ is a bijection between finite quivers on $n$ vertices with no loops and no oriented $2$--cycles and $n\times n$ skew--symmetric integer matrices. We hence write $\myAA(Q)$ for  $\myAA(B(Q))$ and $\TT_n(Q)$ for $\TT_n(B(Q))$. Notice that, in this notation, formula \eqref{eq:MutIntro} of the introduction expresses the mutation of the cluster variable $x_k$ of  the seed $(H,(x_1,\cdots,x_n))$ of the cluster algebra $\myAA(H)$.
%

Every cluster $\mathbf{u}=(u_1,\dotsc, u_n)$ of $\myAA(Q)$ is a free--generating set of the field $\myFF$ and hence $\myFF\simeq\QQ(u_1,\dotsc, u_n)$. In particular, every cluster variable of $\myAA(Q)$ is a rational function in every such cluster. By the famous  Laurent phenomenon proved by Fomin and Zelevinsky in \cite{FZI} such a rational function is actually a Laurent polynomial. We denote by $X_{k;t}^{B;t_0}$ the Laurent expansion in the seed at vertex $t_0$ of $\TT_n(Q)$ whose exchange matrix is $B$ of the $k$--th cluster variable of the seed at vertex $t$ of $\TT_n(Q)$, for $k\in[1,n]$.

A cluster algebra is called of finite type if it has only finitely many cluster variables. In \cite{FZII} it is shown that $\myAA(Q)$ is of finite type if and only if the cluster pattern $\TT_n(Q)$ contains a Dynkin quiver $H$ (i.e. a diagram of type A, D, or E in the Cartan--Killing classification). In this case, as shown in \cite{FZII}, \cite{CC}, \cite{CK1}, the connection with the representation theory of $H$ is much deeper: there is a bijection between the indecomposable $H$--representations and the non--initial cluster variables of $\myAA(H)$. Such bijection is given in terms of projective varieties called quiver Grassmannians. In \cite{CK2} and \cite{DWZII} such bijection is given also for more general quivers but with some restrictions on the involved representations. We will say more about it in the subsequent sections.

\section{Background on quivers with potentials and their representations}\label{Sec:QuivPot}

In this section we recall some facts about the theory of quivers with potentials developed in \cite{DWZ}.

Let $Q=(Q_0,Q_1,t,h)$ be a finite quiver. As usual,  $Q_0$ denotes the set of vertices of Q, $Q_1$ is the set of edges and every edge $a\in Q_1$ is oriented $\xymatrix@1{t(a)\ar^a[r]&h(a)}$.
The theory of quivers with potentials produces a way to ``mutate" the quiver $Q$. More precisely, what is going to change is the set of arrows of Q while the set of vertices is going to remain fixed. In this section we recall how this idea is formalized. 

Let K be a vector space. The vertex span $R=K^{Q_0}$ is defined as the space of K--functions on $Q_0$. There is a distinguished basis  $\{e_i:\,i\in Q_0\}$ of idempotents of R given by $e_i(j)=\delta_{ij}$ (the Kronecker delta) for $i,j\in Q_0$. The arrow span $A=K^{Q_1}$ is the vector space of $K$--functions on the set of arrows. The space $A$ has the following structure of R--bimodule:  $e_ife_j(a)=e_i(h(a))f(a)e_j(t(a))$ for every $a\in Q_1$. We identify the set of arrows $Q_1$ with a  basis of A and for an arrow $a\in Q_1$ we denote with the same symbol $a$ the corresponding element of A. The d--tensor power $A^d=A\otimes\cdots\otimes A$ of A has a structure of R--bimodule as well. Moreover, there is a block decomposition
$$
A^d=\bigoplus_{i,j} A^d_{ij}
$$
where $A^d_{ij}=e_i A^de_j$. The R--bimodule $A^d_{ij}$ is spanned by the elements  $a_1a_2\dotsm a_d$ such that the $a_i\in A$, $h(a_{i+1})=t(a_i)$ for $i\in[1,d-1]$ and $t(a_d)=j$, $h(a_1)=i$. Such elements are called paths of length d from the vertex j to the vertex i. The path algebra is the tensor algebra
$$
R\langle A\rangle:=\bigoplus_{d=0}^\infty A^d
$$
with the convention that  $A^0=R$. For each $i,j\in Q_0$ the R--bimodule $R\langle A\rangle_{i,j}=e_iR\langle A\rangle e_j$ is spanned by the paths from j to i and the union of all the paths form a basis of $R\langle A\rangle$ called the path basis.

For technical reasons it is convenient to consider the completed path algebra
$$
R\langle\langle A\rangle\rangle=\prod_{d=0}^{\infty}A^d
$$
whose elements are possibly infinite linear combinations of paths. Let
$$
\mathrm{m}=\prod_{d=1}^{\infty}A^d
$$
 be the ideal of (linear combinations of) paths of length bigger or equal than one. The algebra $\AAA$ is a topological algebra with respect to the $\mathrm{m}$--adic topology, i.e. a subset $U$ of $\AAA$ is open if and only if for every $x\in U$ there exists $N>0$ such that $x+\mathrm{m}^N\subset U$.

A cyclic path is a path $a_1\dotsm a_d$ such that $t(a_d)=h(a_1)$. We denote by $A^d_{cyc}$ the span of all the cyclic paths in $A^d$. We define the closed subalgebra $\AAA_{cyc}\subseteq\AAA$
by
$$
\AAA_{cyc}=\prod_{d=1}^{\infty}A^d_{cyc}\;.
$$
A potential  $S$ is an element of $\AAA_{cyc}$. Potentials are usually considered up to cyclic equivalences: we say that two potentials $S,S'\in\AAA_{cyc}$ are cyclically equivalent \cite[definition~3.2]{DWZ} if $S-S'$ belongs to the closure of the span of all the elements  of the form $a_1\dotsm a_d-a_2\dotsm a_da_1$ where $a_1\dotsm a_d$ is a cyclic path.

Given a potenital $S$ in $\AAA$, the pair $(A,S)$ is called a quiver with potential (QP for short) provided that $A$ has no loops (i.e. $A_{i,i}=\{0\}$ for every $i\in Q_0$) and no two cyclically equivalent cyclic paths appear in the decomposition of $S$. 

Let $(A,S)$ and $(A',S')$ be two QPs on the same set of vertices $Q_0$. A right--equivalence between $(A,S)$ and $(A',S')$ is an algebra isomorphism $\varphi:\AAA\rightarrow R\langle\langle A'\rangle\rangle$ such that $\varphi|_R=\textrm{id}$ and $\varphi(S)$ is cyclically equivalent to  $S'$. The notion of right--equivalence is very important in dealing with ``mutations" of QPs that we will recall later in section~\ref{Sec:QPMut}.  The direct sum of $(A,S)$ and $(A',S')$ is defined as $(A\oplus A',S+ S')$. Note that this is well-defined since  $R\langle\langle A\rangle\rangle\oplus R\langle\langle A'\rangle\rangle$ embeds canonically in $R\langle\langle A\oplus A'\rangle\rangle$ as a closed subalgebra.

For an element $\xi\in A^\ast$ we consider the cyclic derivative $\partial_\xi$ as the operator $\AAA_{cyc}\rightarrow\AAA$ defined on a cyclic path $a_1\dotsm a_d\in A^d_{cyc}$ by
\begin{equation}\label{Eq:CyclicDerDef}
\partial_\xi(a_1\dotsm a_d)=\sum_{i=1}^d\xi(a_i)a_{i+1}\dotsm a_da_1\dotsm a_{i-1}.
\end{equation}
Given a potential $S$ on $A$, the Jacobian ideal $J(S)$ is the closure of the (two--sided) ideal in $\AAA$  generated by $\{\partial_{\xi}S:\,\xi\in A^\ast\}$. Notice that the closure (in the m--adic topology) of a subset $U\subset\AAA$ is given by $\overline{U}=\bigcap_{N=0}^\infty U+\textrm{m}^N$.
In particular, the closure of an ideal is again an ideal and hence $J(S)$ is a (two--sided) ideal of $\AAA$. The Jacobian algebra is defined as the quotient algebra $\mathcal{P}(A,S)=\AAA/J(S)$.

We notice that if two QPs $(A,S)$ and $(A',S')$ are right--equivalent then the corresponding Jacobian algebras $\mathcal{P}(A,S)$ and $\mathcal{P}(A',S')$ are isomorphic.

We recall the splitting theorem \cite[theorem~4.6]{DWZ}. Let $(A,S)$ be a QP on some set of vertices $Q_0$.  The trivial part $S^{(2)}\in A^2$ of the potential $S$ is, by definition, the homogeneous component of $S$ of degree two. The QP $(A,S)$ is called reduced if $S^{(2)}=0$. Notice that in a reduced QP (A,S), the cyclic part $A^2_{cyc}$ of degree two of $A$ is allowed to be non--zero, even if $S^{(2)}=0$. The trivial and the reduced arrow span of $(A,S)$ are  the R--bimodules given by:
$$
\begin{array}{cc}
A_{triv}=\partial S^{(2)}&A_{red}=A/A_{triv}
\end{array}
$$
where $\partial S^{(2)}$ is the subspace $\{\partial_\xi S^{(2)}:\,\xi\in A^\ast\}\subseteq A$.
The splitting theorem asserts that every QP $(A,S)$ is right--equivalent to the direct sum of a trivial QP $(A_{triv},S_{triv})$ and a reduced QP $(A_{red},S_{red})$. Moreover, the right--equivalence class of both $(A_{red},S_{red})$ and $(A_{triv},S_{triv})$ is determined by the right--equivalence class of $(A,S)$. The QP $(A_{red},S_{red})$ is called the reduced part of $(A,S)$.

\subsection{QP--representations}\label{Sec:QPRep}
A QP--representation is, by definition, a quadruple $(A,S,M,V)$ where $(A,S)$ is a QP, M is a finite dimensional $\mathcal{P}(A,S)$--module and $V=(V_i)_{i\in Q_0}$ is a collection of finite dimensional vector spaces. Sometimes we say that the pair $\mathcal{M}=(M,V)$ is a decorated representation of the QP $(A,S)$. Thus $V$ is a finite dimensional R--bimodule while $M=(M_i)_{i\in Q_0}$ is a finite dimensional representation of the quiver Q whose arrow span is A, which is annihilated by all the cyclic derivatives of the potential S. For every arrow $a\in A$ we denote by $a_M$ the action of a on M. For every vertex $k\in Q_0$ and arrows a and b such that $h(a)=t(b)=k$ there is an element $\partial_{ba}S\in e_{t(a)}\AAA e_{h(b)}$ from the vertex $h(b)$ to the vertex $t(a)$ defined similarly to \eqref{Eq:CyclicDerDef}. Such an element acts on M as a linear map $\gamma_{ba}=(\partial_{ba}(S))_M:M_{h(b)}\rightarrow M_{t(a)}$. This gives rise to a triangle of linear maps
\begin{equation}\label{Eq:TrM}
\xymatrix{
&M_k\ar[dr]^{\beta_M(k)}&\\
M_{in}(k)\ar[ur]^{\alpha_M(k)}&&M_{out}(k)\ar[ll]^{\gamma_M(k)}
}
\end{equation}
where
$$
\begin{array}{cc}
M_{in}(k)=\bigoplus_{a\in Q_1: h(a)=k}M_{t(a)},&M_{out}(k)=\bigoplus_{b\in Q_1:t(b)=k}M_{h(b)}
\end{array}
$$
and 
$$
\begin{array}{ccc}
\alpha_M(k)={\displaystyle\sum_{a\in Q_1: h(a)=k}a_M},&{\displaystyle\beta_M(k)=\sum_{b\in Q_1:t(b)=k}b_M,}&{\displaystyle\gamma_M(k)=\sum_{a,b:\newline h(a)=t(b)=k}\gamma_{ba}}\,.\end{array}
$$
Moreover, the linear maps satisfy the following relations \cite[lemma~10.6]{DWZ}:
\begin{equation}\label{Eq:RelGammak}
\alpha_M(k)\circ \gamma_M(k)=0=\gamma_M(k)\circ \beta_M(k).
\end{equation}
Given two decorated QP--representations $\mathcal{M}=(M,V)$ and $\mathcal{M}'=(M',V')$ of the same QP $(A,S)$, their direct sum is the decorated representation of $(A,S)$ given by $\mathcal{M}\oplus \mathcal{M}'=(M\oplus M',V\oplus V')$.

A QP--representation $\mathcal{M}=(A,S,M,V)$ is called positive if $V=\{0\}$, and negative if $M=0$. The negative simple representation at vertex k is the negative QP--representation  $\mathcal{S}_k^-=\mathcal{S}_k^-(A,S)=(A,S,\{0\},V)$, whose decoration $V$ consists of a  one dimensional vector space at vertex k and zero elsewhere.

A right--equivalence between two QP--representations $\mathcal{M}=(A,S,M,V)$ and $\mathcal{M}'=(A',S',M',V')$ is a triple $(\varphi,\psi,\eta)$ of maps such that: $\varphi:\AAA\rightarrow R\langle\langle A'\rangle\rangle$ is a right--equivalence between $(A,S)$ and $(A',S')$; $\psi:M\rightarrow M'$ is a vector space isomorphism such that $\psi\circ u_M=\varphi(u)_{M'}\circ \psi$;  $\eta:V\rightarrow V'$ is an isomorphism of $R$--bimodules.

Let $(A,S)$ be a QP and let $(A_{red},S_{red})$ be its reduced part. For every trivial QP $(C,T)$ the natural embedding $R\langle\langle A_{red}\rangle\rangle\rightarrow R\langle\langle A_{red}\oplus C\rangle\rangle$ induces an isomorphism of Jacobian algebras $\mathcal{P}(A_{red},S_{red})\rightarrow \mathcal{P}(A_{red}\oplus C,S_{red}+ T)$ \cite[proposition~4.5]{DWZ}. Let $\varphi:R\langle\langle A_{red}\oplus C\rangle\rangle\rightarrow R\langle\langle A\rangle\rangle$ be a right equivalence between $(A_{red}\oplus C,S_{red}\oplus T)$ and $(A,S)$. Given a QP--representation $\mathcal{M}=(A,S,M,V)$, its reduced part is defined as the QP--representation $\mathcal{M}_{red}=(A_{red},S_{red},M',V)$ where $M'=M$ as K--vector space and for $u\in R\langle\langle A_{red}\rangle\rangle$ the action is given by $u_M=\varphi(u)_M$. The right--equivalence class of $\mathcal{M}_{red}$ is determined by that  of $\mathcal{M}$ \cite[proposition~10.5]{DWZ}.

The $\mathbf{g}$--vector of a QP--representation $\mathcal{M}$ is, by definition, the vector $\mathbf{g}_\mathcal{M}=(g_1,\dotsc, g_n)\in \ZZ^n$ ($n=|Q_0|$) whose k--th component is given by
\begin{equation}\label{Def:Gk}
g_k=\textrm{dim ker} \gamma_M(k)-\textrm{dim}M_{k}+\textrm{dim} V_k
\end{equation}
for every $k=1,2,\dotsc, n$ (in the notations of \eqref{Eq:TrM}). In particular it follows that for every two $QP$--representations $\mathcal{M}$ and $\mathcal{M}'$ we have
\begin{equation}\label{Eq:GMopluM'}
\mathbf{g}_{\mathcal{M}\oplus\mathcal{M}'}=\mathbf{g}_{\mathcal{M}}+\mathbf{g}_{\mathcal{M}'}.
\end{equation}
We notice that if $A$ is acyclic, i.e. $\AAA_{cyc}=\{0\}$, then $\gamma_M(k)=0$ and hence the $\mathbf{g}$--vector of a positive QP--representation $M$  equals
$\mathbf{g}_M=-E_A\mathbf{dim}(M)$ where $E_A=(e_{ij})$ is the Euler matrix of $A$ (see e.g. \cite{ASS}) defined by $e_{ii}=1$ and $e_{kj}=-\textrm{dim}A_{jk}$ ($j\neq k$) and $\mathbf{dim}M=(\textrm{dim} M_i)_{i\in Q_0}$.

\subsection{Mutations of QPs}\label{Sec:QPMut}

We  recall the mutation of a quiver with potential $(A,S)$ on some set of vertices $Q_0$. Let $k\in Q_0$ be a vertex such that no oriented 2--cycles in A start (and end) at $k$, i.e. either $A_{ik}=0$ or $A_{ki}=0$ for all $i\in Q_0$. Let us also assume that there are no components of the potential $S$ that start (and end) at the vertex k (if this is the case, it is sufficient to replace $S$ with a cyclically equivalent potential). We define the ``premutation"  of $(A,S)$ as the QP $(\tilde{A},\tilde{S})$  on the same set of vertices $Q_0$ as $(A,S)$ defined as follows: the new arrow span $\tilde{A}$ is given in three steps:
\begin{enumerate}
\item take all the arrows of A which do not start or end at k;
\item for every path $ba$ such that $h(a)=t(b)=k$ add a new arrow $[ba]\in e_{h(b)}\tilde{A}e_{t(a)}$; \item replace every arrow $a$ in $e_kA$ (i.e. ending at k) or in $Ae_k$ (i.e. starting at k) by an opposite arrow $a^\ast$.
\end{enumerate}
The potential $\tilde{S}$ on $\tilde{A}$ is given by
$$
\tilde{S}=\left[S\right]+\Delta_k
$$
where $\left[S\right]$ is obtained by replacing in S every path $ba$ such that $h(a)=t(b)=k$, with the arrow $[ba]$ (recall that there are no components of S starting at $k$); the element $\Delta_k$ is defined by
$$
\Delta_k=\sum a^\ast [ba] b^\ast,
$$
where the sum is taken over all the paths $ba$ such that $h(a)=t(b)=k$. Now the mutation $\mu_k(A,S)$ of the QP $(A,S)$ at vertex $k$ is defined as the reduced part $(\tilde{A}_{red},\tilde{S}_{red})$ of $(\tilde{A},\tilde{S})$. In view of the Splitting Theorem, the operation $(A,S)\mapsto \mu_k(A,S)$ is well--defined on the set of right equivalence classes of QPs.

We remark that by \cite[theorem~5.7]{DWZ} $\mu_k^2(A,S)$ is right--equivalent to $(A,S)$ and hence $\mu_k$ is an involution on the set of right--equivalence classes of QPs.

In order to perform mutations of a QP $(A,S)$ in all the vertices, the arrow span $A$ is assumed to be 2--acyclic i.e. for every vertex $k$ either $A_{ik}=\{0\}$ or $A_{ki}=\{0\}$ for every vertex $i$. Even in this case the mutation can create 2--cycles. A QP is called non--degenerate if this does not happen and any sequence of mutations does not create 2--cycles. If the field K is uncountable then for every arrow span A there exists a potential S such that (A,S) is non--degenerate \cite[corollary~7.4]{DWZ}.

We notice that a 2--acyclic arrow span $A$ with no loops, can be encoded by a $n\times n$ skew--symmetric integer matrix $B=B(A)$ whose ij--th component $b_{ij}$ is given by
$$
b_{ij}=\textrm{dim}A_{ij}-\textrm{dim}A_{ji}.
$$
It can be shown \cite[proposition~7.1]{DWZ} that in this case the mutation $\mu_k(A,S)=(A',S')$ translates into the matrix mutation $B(A')=\mu_k(B(A))$ given by \eqref{Eq:MatrixMut}.

A QP $(A,S)$ is called rigid if every potential $S'$ on A is cyclically equivalent to an element of J(S). Rigid potentials have several nice properties: if (A,S) is rigid then also $\mu_k(A,S)$ is rigid \cite[corollary~6.11]{DWZ}; moreover, rigid potentials are 2--acyclic \cite[proposition~8.1]{DWZ}. It follows that rigid potentials are non--degenerate.

If $Q$ is acyclic, with arrow span A, the only possibility is that  the potential is zero and the QP (A,0) is rigid. Moreover  for every arrow span $A'$ such that $B(A')$ is mutation--equivalent via \eqref{Eq:MatrixMut} to $B(A)$, there exists a potential $S'$ such that $(A',S')$ is reduced and rigid and $(A',S')$ is unique up to right--equivalences. In the special case of a Dynkin quiver $H$ such choice for the potential  $S'$  is explicitly given in \cite[section~9]{DWZII}. Let us recall it.

Let H be a Dynkin quiver and let A be an arrow span on the set of vertices of H such that the corresponding matrix $B=B(A)$ is mutation equivalent to $B(H)$. In \cite{FZII} it is shown that the $ij$--th component $b_{ij}$ of $B$ satisfies $b_{ij}\leq 1$, i.e. the space $A_{ij}$ is either zero or it is one dimensional. For $d\geq3$, a d--chordless cycle in A is a d--cycle whose vertices can be labeled by $\nicefrac{\ZZ}{d\ZZ}$ so that the edges are precisely labeled by pairs $\{i,i+1\}$, $i\in\nicefrac{\ZZ}{d\ZZ}$. In \cite[proposition~9.7]{FZII}, it is shown that all the d--chordless cycles of $A$ are cyclically oriented. A potential S of A is called primitive if it is a linear combination of  all the chordless cycles of A. In \cite[proposition~9.1]{DWZ} it is shown that the QP $(A,S)$ is rigid for every primitive potential S of A. 
In type A and D this is a special case of a general construction due to Labardini--Fragoso \cite{Daniel1}.

\subsection{Mutations of QP--representations}\label{SubSec:MutRep}


Let (A,S) be a non--degenerate QP on a set of vertices $Q_0$ and let  $k\in Q_0$. Let $\mathcal{M}=(M,V)$ be a decorated representation of (A,S). We are going to define a QP representation $\mu_k(\mathcal{M})$ which is a decorated representation of the mutated QP $\mu_k(A,S)$.  We define $\overline{\mathcal{M}}=(\overline{M}=(\overline{M}_i:\,i\in Q_0), \overline{V}=(\overline{V}_i:\, i\in Q_0))$ by
$$
\begin{array}{ccc}
\overline{M}_i=M_i,&\overline{V}_i=V_i&(i\neq k)
\end{array}
$$
and
\begin{equation}\label{eq:Def:MutMk}
\displaystyle{
\begin{array}{cc}
\overline{M}_k=\displaystyle\textrm{im}\gamma\oplus{\textrm{ker }\alpha_M(k)\over \textrm{im } \gamma_M(k)}\oplus {\textrm{ker } \gamma_M(k)\over \textrm{im } \beta_M(k)}\oplus V_k,&\displaystyle \overline{V}_k={
{\textrm{ker } \beta_M(k)\over \textrm{ker } \beta_M(k)\cap \textrm{im } \alpha_M(k)}}.
\end{array}}
\end{equation}
The action of an arrow $c\in \tilde{A}$ (see section~\ref{Sec:QPMut}) on $\overline{M}$ is defined as follows. If c is not incident to k then $c_{\overline{M}}=c_M$ and $[ba]_{\overline{M}}=b_M\circ a_M$ for every arrows a and b of A such that $h(a)=t(b)=k$. It remains to define the linear maps
\begin{equation}\label{Eq:TrMBar}
\xymatrix{\overline{M}_{out}=M_{in}(k)&&\overline{M}_k
\ar_-{\overline{\beta}=(\beta_1,\beta_2,\beta_3,\beta_4)}[ll]&&\overline{M}_{in}=M_{out}(k)\ar_-{\overline{\alpha}=\left(\tiny{\begin{array}{c}\alpha_1\\\alpha_2\\\alpha_3\\\alpha_4\end{array}}\right)}[ll]}
\end{equation}
in the corresponding triangle \eqref{Eq:TrM} for $\overline{\mathcal{M}}$ at vertex k (in \eqref{Eq:TrMBar} we express both $\overline{\alpha}$ and $\overline{\beta}$ in the matrix form corresponding to the decomposition \eqref{eq:Def:MutMk} of $\overline{M}_k$). These maps are defined in the following natural way: we choose splitting data $\rho: M_{out}(k)\rightarrow \textrm{ker }\gamma_M(k)$ such that $\rho\circ \iota=\textrm{id}_{\textrm{ker }\gamma_M(k)}$ (where $\iota$ denotes the inclusion map) and $\sigma:\nicefrac{\textrm{ker }\alpha_M(k)}{\textrm{im }\gamma_M(k)}\rightarrow\textrm{ker }\alpha$ such that $\pi\circ\sigma=\textrm{id}_{\nicefrac{\textrm{ker }\alpha_M(k)}{\textrm{im }\gamma_M(k)}}$ (where $\pi$ denotes the natural projection). The components $\alpha_i$ and $\beta_i$ are defined by:
$$
\begin{array}{ll}
\beta_1=\iota,&\alpha_1=-\gamma\\
\beta_2=\iota\circ\sigma,&\alpha_2=0\\
\beta_3=0,&\alpha_3=-\pi\circ\rho\\
\beta_1=0,&\alpha_1=0
\end{array}
$$
This choice makes $\overline{\mathcal{M}}$ a decorated representation of the premutation $(\tilde{A},\tilde{S})$ of $(A,S)$ \cite[proposition~10.7]{DWZ}. Moreover, a different choice of the splitting data $\rho$ and $\iota$ would produce an isomorphic representation \cite[proposition~10.9]{DWZ}. The mutation of $\mathcal{M}$ is $\mu_k(\mathcal{M}):=\overline{\mathcal{M}}_{red}$, i.e. it is the reduced part of $\overline{\mathcal{M}}$ and hence a decorated representation of $\mu_k(A,S)$. The right--equivalence class of $\mu_k(\mathcal{M})$ is determined by the right--equivalence class of $\mathcal{M}$ \cite[proposition~10.10]{DWZ}. 
Moreover, $\mu_k^2(\mathcal{M})$ is right--equivalent to $\mathcal{M}$ \cite[theorem~10.13]{DWZ}.

\subsection{Some mutation--invariants}\label{SubSec:Invariants}

Let $\mathcal{M}=(M,V)$ and $\mathcal{N}=(N,W)$ be decorated representations of the same nondegenerate QP $(A,S)$. We consider the following number \cite[section~7]{DWZII}
\begin{equation}\label{Eq:DefEInj}
E^{inj}(\mathcal{M},\mathcal{N})=\textrm{dim}\textrm{Hom}_{\mathcal{P}(A,S)}(M,N)+\mathbf{dim}(M)\cdot\mathbf{g}_{\mathcal{N}}
\end{equation}
where $\mathbf{dim}(M)=(\textrm{dim }M_i)_{i\in Q_0}$ is the dimension vector of the positive part M of $\mathcal{M}$, $\mathbf{g}_\mathcal{N}$ is the $\mathbf{g}$--vector of $\mathcal{N}$ whose k--th entry is given by \eqref{Def:Gk} and $\cdot$ denotes the usual scalar product of vectors.

The $E$--invariant of a QP--representation $\mathcal{M}$ is the number
$$
E(\mathcal{M}):=E^{inj}(\mathcal{M},\mathcal{M}).
$$
By \cite[theorem~7.1]{DWZII}, this number is invariant under mutations, i.e. $E(\mu_k(\mathcal{M}))=E(\mathcal{M})$. We notice that $E(\mathcal{S}_k^-,\mathcal{S}_k^-)=0$ and hence $E(\mathcal{M})=0$ for every QP--representation mutation equivalent to $\mathcal{S}_k^-$.

We now recall the homological interpretation of the E--invariant given in \cite[section~10]{DWZII}. Let us assume that the QP $(A,S)$ has the following property:
\begin{equation}\label{CondEinj}
\begin{array}{l}
\text{the potential $S$ belongs to the path algebra $R\langle A\rangle$, and  the two--sided} \\
\text{ideal $J_0$ of $R\langle A\rangle$ generated by all the cyclic derivatives $\partial_aS$}\\\text{contains some power $\textrm{m}^N$ of the ideal m.}
\end{array}
\end{equation}
Under these assumption the Jacobian algebra coincides with $R\langle A\rangle/J_0$ and it is finite dimensional. Condition~\eqref{CondEinj} is satisfied by a rigid QP $(A,S)$ mutation equivalent to a Dynkin quiver in view of  the explicit description of all such QPs recalled above.

Let $\mathcal{M}=(M,0)$ and $\mathcal{N}=(N,0)$ be two positive representations of a QP $(A,S)$ which satisfy \eqref{CondEinj}. The following formula holds \cite[corollary~10.9]{DWZII}
\begin{equation}\label{Eq:EInj}
E^{inj}(M,N)=\textrm{dim}\textrm{Hom}(\tau^{-1}N,M).
\end{equation}
where $\tau$ is the Auslander--Reiten translate.

\section{Quiver with potentials and cluster algebras}\label{Sec:QuivPotClustAlg}

Let K be the field of complex numbers. Let $n\geq1$ and let $\TT_n(B)$ be a cluster pattern as in section~\ref{Sec:CluAlg} associated with a skew--symmetric $n\times n$ integer matrix B. Let A be a set of arrows such that $B(A)=B$ and let S be a potential in $\AAA_{cyc}$ such that $(A,S)$ is a non--degenerate QP. In this section we recall the following two results from \cite{DWZII}:
\begin{enumerate}
\item for every decorated representation $\mathcal{M}$ of $(A,S)$ there is a corresponding a Laurent polynomial $X_{\mathcal{M}}\in \myFF=\QQ(x_1,\dotsc, x_n)$.
\item There exists a family  $\{\mathcal{M}_{k,t}^{B,t_0}\}$ of decorated representations of $(A,S)$ such that
\begin{itemize}
\item for every $k\in[1,n]$, $\mathcal{M}_{k,t_0}^{B,t_0}=\mathcal{S}_k^-=\mathcal{S}_k^-(A,S)$ (defined in section~\ref{Sec:QPRep});
\item for every vertex $t\in\TT_n(B)$ and every index $k\in[1,n]$ we have
$$X_{\mathcal{M}_{k,t}^{B,t_0}}=X_{k,t}^{B,t_0}$$ (defined in section~\ref{Sec:CluAlg}). In particular, $X_{\mathcal{S}_{k}^-}=x_k$.\end{itemize}
\end{enumerate}
Let us start by recalling 1). Let $M$ be a finite--dimensional (complex) representation of a finite quiver $Q$. Given a dimension vector $\mathbf{e}\in\ZZ_{\geq0}^n$, the quiver Grassmannian $Gr_\mathbf{e}(M)$ is the collection of all subrepresentations of $M$ of dimension vector $\mathbf{e}$. It is closed inside the product of the usual Grassmannians $\prod_{i\in Q_0} Gr_{e_i}(M_i)$ and it is hence a projective variety. We denote by $\chi(Gr_\mathbf{e}(M))$ its Euler--Poincar\'e characteristic. The F--polynomial $F_\mathcal{M}$ of a QP--representation $\mathcal{M}=(A,S,M,V)$ is defined as the generating function of $\chi(Gr_\mathbf{e}(M))$:
\begin{equation}\label{Def:FM}
F_M(y_1,\dotsc, y_n)=\sum_{\mathbf{e}}\chi(Gr_\mathbf{e}(M))y_1^{e_1}\dotsm y_n^{e_n}.
\end{equation}
In particular, the $F$--polynomial of a negative QP--representation is 1. In 
\cite[proposition~3.2]{DWZII}, it is shown that
\begin{equation}\label{FMoplusM'}
F_{\mathcal{M}\oplus \mathcal{M}'}=F_\mathcal{M}F_{\mathcal{M}'}.
\end{equation}
Let $\mathbf{b}_1,\dotsc,\mathbf{b}_n$ be the columns of the matrix B.
The desired Laurent polynomial $X_{\mathcal{M}}$ is defined by:
\begin{equation}\label{Eq:DefXM}
X_{\mathcal{M}}=F_{\mathcal{M}}(\mathbf{x}^{\mathbf{b}_1},\dotsc, \mathbf{x}^{\mathbf{b}_n})\mathbf{x}^{\mathbf{g}_\mathcal{M}},
\end{equation}
where $\mathbf{g}_\mathcal{M}$ is the $\mathbf{g}$--vector of $\mathcal{M}$ defined in \eqref{Def:Gk} and we use the notation $\mathbf{x}^{\mathbf{c}}=x_1^{c_1}\dotsm x_n^{c_n}$ for $\mathbf{c}=(c_1,\dotsc,c_n)$.
In particular, we can rewrite \eqref{Eq:DefXM} as follows:
\begin{equation}\label{Eq:XM}
X_\mathcal{M}=\sum_{\mathbf{e}} \chi(Gr_{\mathbf e}(M)) {\mathbf x}^{ \mathbf g_{M}+B\mathbf{e}}.
\end{equation}
From \eqref{FMoplusM'} and \eqref{Eq:GMopluM'} it follows that  the map $\mathcal{M}\mapsto X_{\mathcal{M}}$ has the following property:
\begin{equation}\label{Eq:XMoplusM'}
X_{\mathcal{M}\oplus \mathcal{M}'}=X_\mathcal{M}X_{\mathcal{M}'}
\end{equation}
for every decorated representation $\mathcal{M}$ and $\mathcal{M}'$ of $(A,S)$. Moreover, it can be shown that the map $\mathcal{M}\mapsto X_\mathcal{M}$ is constant on right--equivalence classes of QP--representations.

The expression \eqref{Eq:XM} is a sum of Laurent monomials and hence it is a Laurent polynomial in the variables $x_1,\dotsc,x_n$. Its reduced form is hence a rational function whose denominator is a monomial $x_1^{d_1}\dotsm x_n^{d_n}$. The integer vector $\mathbf{d}(X_\mathcal{M})=(d_1,\dotsc, d_n)$ is called the denominator vector of $X_\mathcal{M}$. For example, $\mathbf{d}(X_{\mathcal{S}_k^-})=\mathbf{d}(\nicefrac{1}{x_k^{-1}})=(0,\dotsc,-1,\dotsc,0)$, -1 at the k--th position. In \cite[corollary~5.5]{DWZII}, it is shown that
\begin{equation}\label{Eq:DenDim}
d_i\leq dim (M_i)\quad \textrm{for every } i\in[1,n].
\end{equation}

Let us recall 2). Let $t$ be a vertex of $\TT_n(B)$ and let $\Sigma_t=(B',\mathbf{u})$ be the corresponding seed (see section~\ref{Sec:CluAlg}). Since $\TT_n(\Sigma)$ is a tree, there is a unique path $$\xymatrix@1{t_0\ar@{-}^{k_1}[r]&t_1\ar@{-}^{k_2}[r]&\cdots\ar@{-}^{k_{m-1}}[r]&t_{m-1}\ar@{-}^{k_{m}}[r]&t_m=t}$$ which connects $t_0$ with $t$. Let
\begin{equation}\label{Eq:DefAtSt}
(A_t,S_t)=\mu_{k_m}\circ\cdots\circ\mu_{k_2}\circ\mu_{k_1} (A,S).
\end{equation}
In particular $(A_t,S_t)$ is a non--degenerate QP and $B(A_t)=B'$.
Recall that  $X_{k,t}^{B,t_0}$ denotes the Laurent expansion of the k--th cluster variable $u_k$ of $\Sigma_t$ in the seed $\Sigma_{t_0}=(B,\mathbf{x})$.  We have
$$
u_k=X_{k,t}^{B,t_0}=\mu_{k_m}\circ\cdots\circ\mu_{k_2}\circ\mu_{k_1} x_k
$$
Let $\mathcal{S}_k^-=\mathcal{S}_k^-(A_t,S_t)$ be the negative simple representation of $(A_t,S_t)$. We define
\begin{equation}\label{Eq:DefMutCl}
\mathcal{M}_{k,t}^{B,t_0}:=\mu_{k_1}\circ\cdots\circ\mu_{k_{m-1}}\circ\mu_{k_m} \mathcal{S}_k^-
\end{equation}
By \cite[theorem~5.1]{DWZII}, $X_{\mathcal{M}_{k,t}^{B,t_0}}=X_{k,t}^{B,t_0}$. The family $\{\mathcal{M}_{k,t}^{B,t_0}\}$ is the desired family.

We notice that by \eqref{Eq:XMoplusM'} the same description holds for cluster monomials. For instance the Laurent expansion of a cluster monomial $b=u_1^{a_1}\dotsm u_n^{a_n}$ is given by
\begin{equation}\label{Eq:ClMonExp}
b=X_{\mathcal{M}_{1,t}^{B,t_0\oplus a_1}\oplus\mathcal{M}_{2,t}^{B,t_0\oplus a_2}\oplus\cdots\oplus \mathcal{M}_{n,t}^{B,t_0\oplus a_n}}
\end{equation}
and we define
\begin{equation}\label{Eq:ClMonRep}
\mathcal{M}_b^{B,t_0}=\mathcal{M}_{1,t}^{B,t_0\oplus a_1}\oplus\mathcal{M}_{2,t}^{B,t_0\oplus a_2}\oplus\cdots\oplus \mathcal{M}_{n,t}^{B,t_0\oplus a_n}.
\end{equation}
The representations $\mathcal{M}_{k,t}^{B,t_0}$ have the following remarkable property: in view of \eqref{Eq:DefMutCl}, $E(\mathcal{M}_{k,t}^{B,t_0})=0$ for every $k$ and $t$.  For a QP--representation $\mathcal{M}$ such that $E(\mathcal{M})=0$ we have \cite[corollary~5.5]{DWZII}:
\begin{equation}\label{Eq:E0Cor}
\textrm{either }M_i=\{0\}\textrm{ or }V_i=\{0\}\textrm{ for every }i\in[1,n].
\end{equation}
In particular, this holds for the family $\{\mathcal{M}_{k,t}^{B,t_0}\}$.

It is remarkable that, while  the definition of the family $\{\mathcal{M}_{k,t}^{B,t_0}\}$ depends on the choice of a non--degenerate QP (A,S), the cluster algebra $\myAA(B)$ only depends on the initial exchange matrix $B=B(A)$. In general, two potentials $S$ and $S'$ on A such that  $(A,S)$ and $(A,S')$ are non--degenerate can be very different and give rise to non--isomorphic Jacobian algebras. 

\section{Proof of theorem~\ref{MainThmIntro}}\label{Sec:Proof}
In view of \cite[corollary~3]{CK1}, the set $\BB$ of cluster monomials form a $\ZZ$--basis of $\myAA(H)$. In view of \cite{MS, hernandez, Nakajima}  the cluster monomials are positive. We hence prove that $\BB$ is the atomic basis of $\myAA(H)$ i.e.  positive elements of $\myAA(H)$ are non--negative linear combinations of cluster monomials.
We say that a Laurent monomial $x_1^{a_1}\dotsm x_n^{a_n}$ in some variables $x_1,\dotsc, x_n$ is \emph{proper} if there is at least an index $i$ such that $a_i<0$.
A shown in \cite{Sherman} the following lemma implies theorem~\ref{MainThmIntro}. 

\begin{lemma}\label{lemma}
For every cluster $\myCC$ of $\myAA(H)$ and every cluster monomial $b$ which is not a monomial in the elements of $\myCC$, the expansion of $b$ in $\myCC$ is a sum of proper Laurent monomials.
\end{lemma}
We hence prove lemma~\ref{lemma}. Let $\TT_n(H)$ be the cluster pattern associated with H (see section~\ref{Sec:CluAlg}). In particular, n denotes the number of vertices of H. Let $\Sigma_0=(H,\mathbf{x})$ be the initial seed at some vertex s of $\TT_n(H)$. Let $t_0$ be a vertex of $\TT_n(H)$ and let $\Sigma_{t_0}=(B,\mathbf{u})$ be the corresponding seed in $\TT_n(H)$. To such a vertex, there is also associated a QP $(A,S)=(A_{t_0},S_{t_0})$, which is mutation equivalent to the QP $(H,0)$ by \eqref{Eq:DefAtSt}. In section~\ref{Sec:QPMut} we have recalled the explicit description of  $(A,S)$ and we have noticed that it is rigid and hence non--degenerate.

Let $b$ be a cluster monomial of $\myAA(H)$ and let $\mathcal{M}_b^{B,t_0}$ be the corresponding decorated representation of $(A,S)$ given by \eqref{Eq:ClMonRep}. We show that the Laurent polynomial $X_{\mathcal{M}_b}^{B,t_0}=b$ given by \eqref{Eq:ClMonExp} is a sum of proper Laurent monomials.

Let us first consider the case when $\mathcal{M}_b^{B,t_0}=(M,V)$ is not a positive representation, i.e. there is an $i\in Q_0$ such that $V_i\neq\{0\}$. Therefore the monomial $b$ has the form $x_ib'$ for another cluster monomial $b'$. For such an index i, since $E(\mathcal{M}_b^{B,t_0})=0$ and in view of \eqref{Eq:E0Cor}, we have that $M_i=\{0\}$. In view of \eqref{Eq:DenDim}, the i--th entry $d_i$ of the denominator vector of $X_{\mathcal{M}_b^{B,t_0}}$  is zero as well. It follows that if we prove the lemma for $b'$ then the lemma holds also for $b$.

We hence assume that $\mathcal{M}_b^{B,t_0}$ is a positive representation of $(A,S)$, i.e. $\mathcal{M}_b^{B,t_0}=(M,0)$ and $M$ is a finite--dimensional $\mathcal{P}(A,S)$--module. In view of \eqref{Eq:XM}, we have
$$
X_{\mathcal{M}_b^{B,t_0}}=\sum_{\mathbf{e}} \chi(Gr_{\mathbf e}(M)) {\mathbf x}^{ \mathbf g_{M}+B\mathbf{e}}.
$$
We show that if there exists a non--zero subrepresentation $N$ of $M$ of dimension vector $\mathbf e$ then the vector $\mathbf g_{M}+B\mathbf{e}$ has at least one negative entry. Since $B$ is skew--symmetric, the scalar product $\mathbf e\cdot(\mathbf g_{M}+B\mathbf{e})=\mathbf e\cdot\mathbf{g}_M$. We hence show that the number $\mathbf e\cdot\mathbf{g}_M$ is negative. In view of \eqref{Eq:EInj}  we have
$$
E(M):=\mathbf{dim}(M)\cdot\mathbf g_M+\textrm{dim}\textrm{Hom}(M,M)=0=\textrm{dim}\textrm{Hom}(\tau^{-1}M,M).
$$
Moreover, again by \eqref{Eq:EInj}, we have
$$
E^{inj}(N,M):=\mathbf e\cdot\mathbf g_M+\textrm{dim}\textrm{Hom}(N,M)=\textrm{dim}\textrm{Hom}(\tau^{-1}M,N)
$$
Since $N$ is a subrepresentation of $M$, there is an injection
$$
\textrm{\textrm{Hom}}(\tau^{-1}M,N)\rightarrow\textrm{\textrm{Hom}}(\tau^{-1}M,M)
$$
and so $E^{inj}(N,M)\leq E(M)=0$. It follows that $E^{inj}(N,M)=0$ and hence $\mathbf e\cdot\mathbf g_M=-\textrm{dim}\textrm{Hom}(N,M)<0$ as desired.
It remains to discuss the case $\mathbf{e}=0$. We hence prove that the vector $\mathbf{g}_{\mathcal{M}}$ has at least one negative entry. We take the scalar product $\mathbf{dim}(M)\cdot \mathbf{g}_{\mathcal{M}}=-\textrm{dim}\textrm{Hom}(M,M)<0$, as desired.

\begin{remark}\label{Rem:Coeff}
In view of the separation formula \cite[corollary~6.3]{FZIV} and of the explicit expression for the F--polynomials and the $\mathbf{g}$--vectors of a cluster monomial in every cluster (recalled in section~\ref{Sec:QuivPotClustAlg}) we can prove lemma~\ref{lemma} in a cluster algebra $\myAA_\PP(H)$ with arbitrary coefficients $\PP$. We know that cluster monomials of $\myAA_\PP(H)$ are positive (i.e. their Laurent expansions in every cluster have coefficients in $\ZZ_{\geq0}\PP$). If we assume that the they form a $\ZZ\PP$--basis of $\myAA_\PP(H)$ then they are an atomic basis for the same reasons as in the coefficient--free setting (see \cite{Sherman}). 
\end{remark}

\bibliographystyle{plain}
\bibliography{Bibliografia10Atom}

\end{document}